\documentstyle[12pt,amstex,psfig]{article}

\newtheorem{prop}{Proposition}
\newtheorem{lemma}[prop]{Lemma}
\newtheorem{cor}[prop]{Corollary}
\newtheorem{thm}[prop]{Theorem}
\newtheorem{oldthm}[prop]{Theorem}
\newtheorem{prob}{Problem}
\newcommand{\R}{{\Bbb R}}

\newcommand{\Z}{{\Bbb Z}}

\begin{document}
\title{Tiling a rectangle with the fewest squares.}
\author{Richard Kenyon\thanks{CNRS UMR 128, Ecole Normale Sup\'erieure de Lyon,
46, all\'ee d'Italie, 69364 Lyon, France. Research at MSRI supported in part
by NSF grant no. DMS-9022140.}}
\date{}
\maketitle
\begin{abstract}
We show that a square-tiling of a $p\times q$ rectangle, 
where $p$ and $q$ are relatively prime integers,
has at least $\log_2p$ squares. If $q>p$ we construct a square-tiling
with less than $q/p+C\log p$ squares of integer size, 
for some universal constant $C$.
\end{abstract}

\section{Introduction}
A certain store sells square tiles of arbitrary positive integer size
for \$1 each. You'd like to tile your kitchen
(a $p\times q$ rectangle, $p,q\in\Z_+$), for the least cost. What's the cheapest way?

We show:
\begin{thm} A $p\times q$ rectangle, where $p,q$ are relatively prime integers, $p<q$,
requires at least $\max\{q/p,\log_2q\}$ square tiles to tile.
Furthermore there exists a square tiling with less than $q/p+C_1\log_2p$ squares
of integer size, for some universal constant $C_1$. 
\label{int}
\end{thm}

\noindent{\bf Remark.} 
For the lower bound the sizes of the squares are not restricted to be integers.
Also, the quantity $q/p$ in the two bounds is
necessary for thin rectangles;
for example an $n\times 1$ rectangle requires at least $n$ squares. 
If $q/p$ is bounded then we have logarithmic upper and lower
bounds.
\medskip

Here for a $p\times q$ rectangle we call the {\bf aspect ratio} the larger of $p/q,q/p$. 

In case the aspect ratio $x>1$ of the kitchen is not rational, no tiling with a finite number of squares 
is possible by a theorem of Dehn \cite{Dehn}; on the other hand, using squares of
arbitrary real size, if you have a 
refrigerator to cover up the untiled portion, you can do equally well:
\begin{thm} 
For any $\epsilon>0$ and $x\in\R,~x>1$, one can tile all but an $\epsilon$-
neighborhood of a corner of an $x\times 1$ rectangle with 
$\leq x+C_2\log(1/\epsilon)$ squares, for some universal constant $C_2$.
\label{2}
\end{thm}

Our proof of the lower bound in Theorem \ref{int} uses the theory of electrical
networks, which has a well-known connection with square tilings \cite{Dehn,BSST,Bol}. 
In particular a generalization of Theorem \ref{int} is as follows:
\begin{thm} 
Let $X$ be a resistor network with underlying graph $G$, with resistances $1$ on each edge.
If the effective
resistance between two vertices is the rational number $\frac qp>1$
in lowest terms,
then there are at least $\max\{q/p,\log_2q\}$ edges in $G$. Conversely,
for any rational $\frac qp>1$ there is such 
a network with a {\bf planar} graph
having at most $\frac qp+C_1\log_2p$ edges.
\end{thm}

Here we allow multiple edges between the same vertices in the graph $G$.

\section{The greedy algorithm}
Your initial reaction is of course to tile using the greedy algorithm,
that is, select the largest square that fits (a $p\times p$ tile), place it touching
a shortest side of the kitchen, and repeat with the remaining untiled part, which is
now $p\times (q-p)$.

\begin{figure}[htbp]
\centerline{\psfig{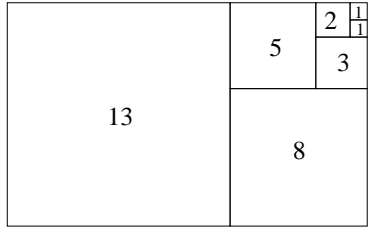}}
\caption{\label{fib} Tiling Fibonacci's kitchen.}
\end{figure}

This method, also known as the {\bf Euclidean algorithm}, works well for certain
shapes of rectangle, for example those rectangles which are $F_n\times F_{n+1}$, where
$F_n$ is the $n$th Fibonacci number. Indeed such a rectangle is tiled with $n\approx
\log_\tau F_n$ tiles (Figure \ref{fib}).

Unfortunately for many shapes of rectangles, this algorithm is quite expensive:
for a $p\times (p+1)$ rectangle, the first square leaves a $p\times 1$ rectangle,
which requires at least $p$ squares to tile, for a total cost
of $1+p$, which much more expensive than is necessary.

We leave the reader to verify that, if the continued fraction
expansion of $q/p$ is
$[a_0;a_1,\ldots,a_k]$, that is,
$$\frac qp=a_0+\frac{1}{a_1+\frac{1}{a_2+\ldots+\frac{1}{a_k}}},$$
then the cost of the greedy algorithm on a $p\times q$ rectangle
is $a_0+a_1+\ldots +a_k$.
(Here the integers $a_i$ are called {\bf partial quotients} of $q/p$: it is required that
$a_0\geq 0$ and for $i\geq 1$ that $a_i\geq 1$. Under these conditions
the $a_i$ are uniquely defined
except for $a_k$, and we have $[a_0;a_1,\ldots,a_k]=[a_0;a_1,\ldots,a_k-1,1]$ assuming
$a_k>1$.)
\medskip

For the irrational rectangle $x\times 1$,
if $x=[a_0;a_1,a_2,\ldots]$ is the infinite continued fraction expansion of $x$, 
then the cost $T_\epsilon(x)$ to cover up all but the $\epsilon$-neighborhood
of the corner is 
\begin{equation}T_\epsilon(x)=a_0+a_1+\ldots+a_k,
\label{cex}
\end{equation}
where the $k$ is the first number to satisfy
$$|x-\frac{p_k}{q_k}|<\frac{\epsilon}{q_k},$$
where $p_k/q_k=[a_0;a_1,a_2,\ldots,a_k]$ is the $k$th rational approximant to $x$.
Indeed, if the untiled portion after $a_0+\ldots+a_k$ steps is $x'\times y'$, then we have
$$\left(\begin{array}{cc}1&-1\\0&1\end{array}\right)^{a_k}\ldots
\left(\begin{array}{cc}1&0\\-1&1\end{array}\right)^{a_1}
\left(\begin{array}{cc}1&-1\\0&1\end{array}\right)^{a_0}
\left(\begin{array}{c}x\\1\end{array}\right)=
\left(\begin{array}{c}x'\\y'\end{array}\right),$$
or in other words
$$\left(\begin{array}{cc}p_k&p_{k\pm1}\\q_k&q_{k\pm1}\end{array}\right)
\left(\begin{array}{c}x'\\y'\end{array}\right)=
\left(\begin{array}{c}x\\1\end{array}\right).$$
(Here the $\pm$ depends on the parity of $k$.)
If $y'<\epsilon$ then 
$$\epsilon>y'=-q_kx+p_k.$$

Quantities related to this cost $T_\epsilon(x)$ for ``typical'' numbers
have been studied in detail. Yuval Peres combined some known results to prove:
\begin{thm}[Peres]
\label{Per}
There is a constant $c'>0$ such that 
for any $\delta>0$ the Lebesgue measure of the set 
$$\left\{x\in(0,1) :
\left|\frac{T_\epsilon(x)}{\log(1/\epsilon)\log\log(1/\epsilon)}-c'
\right|\geq\delta\right\}$$
tends to zero with $\epsilon$.
\end{thm}

For the proof, see the appendix. 
\medskip

If $x=[a_0;a_1,a_2,\ldots]$ is irrational
and all the $a_i$ are bounded by $n$, then $x$
is called {\bf $n$-aloof}. It is not hard to see (we'll see later in any case) that
if $x$ is $n$-aloof, the greedy algorithm gives a logarithmic bound
$T_\epsilon(x)<const\log(\frac 1\epsilon)$, where the constant depends on $n$.

\section{The lower bound.}
We give here a proof of the lower bound in Theorem \ref{int}.
First, the largest square which can fit is a $p\times p$ square,
which covers $p/q$ of the area, so you need at least $q/p$ squares to tile.
We will show that you need at least $\log_2q$ squares to tile.

To a square tiling of a rectangle $R$, associate a graph $G$ as follows \cite{BSST}:
let $G=(V,E)$ be the graph with vertex set $V$ and edges $E$, where $V$
is the set of connected components of the union of the horizontal boundaries
of tiles in the square-tiling, and $E$ is the set of tiles (note that a tile
connects exactly two horizontal components, and that multiple
edges between two vertices are possible).
The vertex corresponding to the upper boundary of $R$
is called $a$, and the vertex corresponding to the lower boundary is $b$.
It is clear that $G$ is planar, and that
$a$ and $b$ are on the same face (the outer face) of $G$.

It is helpful to direct the edges from the upper vertex to the lower vertex.

Associated to $G$ is the resistor network, obtained by assigning each
edge of $G$ a resistance $1$. By assigning {\bf potentials} $p_a,p_b\in\R$ to the vertices
$a$ and $b$ of the network, a flow of electric current is set up in $G$, that is,
we have maps $p\colon V\to\R$ (``potentials'') and $c\colon E\to \R$  (``currents'')
which satisfy Kirchoff's rule and Ohm's law: 
the net current flow out of any vertex (except $a$ and $b$) is equal to the net flow into 
that vertex,
and the current across an edge equals the drop in potential between its endpoints.
(By definition the current has a sign which depends on the direction of the edge.)
The potentials and currents are the unique solution to the equations arising
from Kirchoff's and Ohm's rules with the given boundary conditions $p(a)=p_a,
p(b)=p_b$.
If we denote by $c_b$ the net current going into $b$, then the quantity
$r(a,b)=(p_a-p_b)/c_b$ depends only on the graph and is independent of $p_a,p_b$. 
This quantity $r(a,b)$ is called the {\bf effective resistance}, or {\bf impedance}, from 
$a$ to $b$.

If we scale the square tiling by a homothety of $\R^2$ 
and translate it so that the upper boundary is at
$y$-coordinate $p_a$ and the lower boundary is at $y$-coordinate $p_b$ (assuming with loss
of generality that $p_a>p_b$), 
then we see that a solution (hence the unique solution)
of Kirchoff's and Ohm's equations is given by:
for $v\in V$, $p(v)$ 
equals the $y$-coordinate of the horizontal component corresponding
to $v$, and for $e\in E$, $c(e)$ is the size of the square tile
corresponding to $e$.
The quantity $r(a,b)$ is simply the ratio of height to width of the rectangle.

This construction also works in the other direction. 
We can associate, to any planar graph $G=(V,E)$ 
and choice of two vertices $a,b\in V$
on the same face, a square-tiling of a rectangle $R$ whose resistor network
is $G$. This is proved in \cite{BSST}, who simply use the idea of the
previous paragraph to construct the tiling from the potentials and currents.

By a result of Kirchhoff \cite{Kirch} (see also \cite{BSST}), the resistance
$r(a,b)$ satisfies: $r(a,b)=\kappa_{ab}/\kappa$, where $\kappa$ 
is the number of spanning trees in $G$, and $\kappa_{ab}$
is the number of spanning trees in $G_{ab}$, the graph obtained from $G$
by gluing together vertices $a$ and $b$. 

Thus if $R$ is a $p\times q$ rectangle, 
we have $p/q=\kappa_{ab}/\kappa$, and since $p$ and $q$ are relatively prime, $\kappa_{ab}\geq p$ and
$\kappa\geq q$.

However the number of spanning trees in any graph of $m$ edges is less than $2^m$, 
since a tree is a subset of edges, and there are 
$2^m$ distinct subsets of $m$ edges.  
Since $G$ has $m$ edges, where $m$ is the number of tiles,
we have $q\leq\kappa\leq 2^m$.
We conclude that $m\geq \log_2q$.

This gives the lower bound in Theorem \ref{int}.

\section{The upper bound for real rectangles.}
Let $R$ be a rectangle with aspect ratio $x$ (recall $x\geq 1$).
We assume $x<2$: if not, apply the greedy algorithm $\lfloor x\rfloor -1$
times. The remaining untiled portion has aspect ratio $x-(\lfloor x\rfloor -1)$
in the range $[1,2)$.

We show how to tile $R$ quickly. Assume $R$ is $x\times 1$.

The idea is simple: use the greedy algorithm, 
getting a nested decreasing sequence of rectangles $R_j$
(the untiled portions), each containing a fixed corner of $R$,
of aspect ratios $x=x_0,x_1,x_2,\ldots$, with 
$$x_{j+1}=\max\{x_j-1,\frac{1}{x_j-1}\}$$
until some $x_i$ is close to $1$,
say $x_i<1+\delta$ for some small $\delta>0$. Note then that 
for $0\leq j\leq i$ we have $x_j<1/\delta$. 

At step $i$, instead of putting in a square, which would result in
the new rectangle having aspect ratio $x_{i+1}>1/\delta$, 
just put in a rectangle of aspect ratio $2$,
with its longer side covering the shorter side of $R_i$. The remaining untiled 
portion is a rectangle $R_{i+1}$ with 
aspect ratio $x_{i+1}=1/(x_i-\frac 12)$, and so $1\leq x_{i+1}<2$. Now continue.

Since for each $j$ we have $1<x_j<\frac 1\delta$, 
each square added removes either a fraction
at least $\delta$ of the area (in the case 
when one square is added to a rectangle of aspect ratio
in $[1+\delta,1/\delta]$), or a fraction of at least
$\frac{1}{4(1+\delta)}$ of the 
remaining area (in case a rectangle of aspect ratio $2$, which is tiled by two
squares, is added to a rectangle of aspect ratio $<1+\delta$). 
So the area decreases by a factor of at least
$$\max\{1-\delta,1-\frac{1}{4(1+\delta)}\}$$ per square added.

This quantity is minimized 
when $\delta=\frac{\sqrt{2}-1}2\approx.207$, and the rate is $\lambda=1-\delta\approx.793$.

After $k$ squares, the untiled area is a rectangle of area at most 
$\lambda^kx$, and aspect
ratio between $1$ and $\frac 1\delta$, 
and so is contained in a neighborhood of 
radius $x\lambda^{k/2}/\delta$ of the corner.

The completes the proof of Theorem \ref{2}.\hfill{$\Box$}
\medskip

\noindent{\bf Remark.} The rate of decrease of area
$\lambda=.793$ is of course not the optimal one.
Optimization of similar algorithms seems to be an interesting problem, but
one we won't consider here. One might conjecture that
$\tau^2=\frac{3-\sqrt{5}}{2}\approx .381$ is a lower
bound for $\lambda$, since the golden rectangle seems to be most easily 
tiled by the greedy algorithm, which has this rate.
\medskip

An alternative method for tiling an $x\times 1$ rectangle 
is suggested by the following result of Hall:
\begin{oldthm}[Hall \cite{Hall}]
Any real number between $\sqrt{2}-1$ and $4+4\sqrt{2}$ 
can be written as the sum of two $4$-aloof numbers.
\label{thick}
\end{oldthm}

Here $\sqrt{2}-1=2\cdot[0;4,1,4,1,\ldots]$ and $4+4\sqrt{2}=2\cdot[4;1,4,1\ldots]$ are
the minimal and maximal possible sums.

So to tile a $x\times 1$ rectangle, where $1< x\leq 4+4\sqrt{2}$,
write $x=x_1+x_2$ where $x_1,x_2$ are $4$-aloof, and divide the $x\times 1$
rectangle into an $x_1\times 1$ and an $x_2\times 1$ rectangle with a single vertical line.
Now tile each of the subrectangles using the greedy method; by $4$-aloofness,
in each subrectangle each new square added takes up
at least $1/5$ of the area, so that after $2n$ squares the remaining area is
at most $(4/5)^n$ of the original area there. The rate is then $\lambda=(4/5)^{1/2}$.

Our method for integral rectangles will be a variant on this method.

\section{An upper bound for rational rectangles.} 
Let $R$ be a $p\times q$ rectangle, with $p,q\in \Z$, $(p,q)=1$, $0<p\leq q$. 
We assume as before that $q/p<2$: if it is larger, use the greedy algorithm
$n=\lfloor q/p\rfloor-1$ times, so that the remaining $p\times (q-np)$ rectangle
satisfies the above conditions.

We establish in this section an upper bound of $C\log p\log\log p$ for the number of
squares needed to tile $R$. Section \ref{log} refines the construction 
to improve the bound to $C_1 \log p$.

The construction proceeds as follows. Let $x_1,x_2$ be $4$-aloof numbers
such that $x_1+x_2=q/p$ (using Theorem \ref{thick}). 
Let $k_1=\lfloor x_1p\rfloor$, and $k_2=q-k_1$,
so that $|\frac{k_i}{p}
-x_i|<\frac 1p$ and $$\frac{k_1}p+\frac{k_2}p=\frac qp.$$

We divide the rectangle $R$ into a $p\times k_1$ rectangle
$R_1$ and a $p\times k_2$
rectangle $R_2$. We will show (below, after Lemma \ref{sqrt})
that we can apply the greedy algorithm successfully to each of these
rectangles for a while, that is,
until the remaining untiled rectangles $R_1',R_2'$ each have 
side lengths $\leq c\sqrt{p}$ for some universal constant $c$ 
(and $R_1',R_2'$ have aspect ratios $\leq 2$).

We then repeat the process, using Theorem \ref{thick} again to subdivide 
$R_1',R_2'$ each in two, applying the greedy algorithm until the remainders
have sides $\leq c\sqrt{cp^{1/2}}\leq c^2p^{1/4}$, and so on.

We show that 
at each stage the number of squares added in a single rectangle before we 
subdivide it is at most the logarithm to the base $\alpha=6/5$ of
its larger side length. The side lengths decrease by at least $x\mapsto c\sqrt{x}$
before we resubdivide, and each subdivision doubles the number of rectangles.
When the edge lengths of a subrectangle are less than the constant $2c^2$ in length,
simply tile the subrectangle in any way you please.

We derive for the total number $N$ of squares needed to tile:
$$N\leq 2\log_\alpha p+4\log_\alpha (cp^{1/2})+8\log_\alpha (c^{3/2}p^{1/4})+
\ldots+2^k\log_\alpha (c^{2-2^{-k+1}}p^{2^{-k}})+2^{k+1}c',$$
where $k$ is chosen so that $p^{2^{-k}}\approx 2,$ that is, $k\approx\log_2\log_2 p$,
and $c'$ is the number of squares needed to tile an integer-sided
rectangle whose sidelengths are bounded by $2c^2$.
(Note that $c^bp^\beta\to c^{1+b/2} p^{\beta/2}$ under the map $x\mapsto c\sqrt{x}$.)

Thus the number of squares is bounded above by
$$N\leq 2k\log_\alpha (c^2p)+2^{k+1}c'\leq2\log_\alpha (c^2p)\log_2\log_2p+2c'\log_2p.$$
\medskip

It remains to prove our claim that we can tile a $p\times k_1$
rectangle quickly using the greedy algorithm until the remining untiled rectangle
has edges $\leq c\sqrt{p}$.

Recall that a {\bf Farey interval} $I$ is a subinterval of $(0,\infty)$
with rational endpoints
$(\frac{p_1}{q_1},\frac{p_2}{q_2})$ which satisfy $p_2q_1-p_1q_2=1$.
(Notationally we allow $p_2/q_2=\infty=``\frac 10"$ and $p_1/q_1=0=``\frac01"$.)
Each Farey interval $I$ gives rise to two Farey subintervals 
$L(I)=(\frac{p_1}{q_1},\frac{p_1+p_2}{q_1+q_2})$ and 
$R(I)=(\frac{p_1+p_2}{q_1+q_2},\frac{p_2}{q_2})$,
and the set of all Farey intervals form a binary tree in this way with the root
being $I_0=(0,\infty)=(\frac{0}{1},\frac{1}{0})$.
A Farey interval has a label indicating the unique descending path to $I$
from the root; this label is a finite word in the letters `L' and `R'.
thus for example $LRL(I_0)=LR((\frac 01,\frac 11))=
L((\frac 12,\frac 11))=(\frac 12,\frac 23)$.
The Farey interval $I=R^{a_0}L^{a_1}R^{a_2}L^{a_3}\ldots R^{a_{2k}}(I_0)$
has the property that for $x\in I$, the continued fraction expansion of $x$
begins $x=[a_0;a_1,a_2,\ldots,a_{2k},\ldots]$, and similarly for
words ending in $L^{a_{2k+1}}$.
We call a Farey interval {\bf finite} if $q_2>0$.

\begin{lemma} If $I=(\frac{p_1}{q_1},\frac{p_2}{q_2})$ 
is a finite Farey interval and contains a $4$-aloof number, 
then $1/5\leq \frac{q_1}{q_2}\leq 5$.
\label{bound}
\end{lemma}
\noindent{\bf Proof.} 
If $x\in I$ is $4$-aloof, the word $w$ such that $I=w(I_0)$
has no more than $4$ consecutive
$L$'s or $R$'s. In particular $1/5<x<5$, so $x$ is in one of the intervals
$$(\frac15,\frac14),(\frac14,\frac13),(\frac13,\frac12),(\frac12,\frac 11),(\frac
11,\frac21),(\frac21,\frac31),(\frac31,\frac41),(\frac41,\frac51),$$
for which the result is true.

Now if $I=L(J)$, then clearly $q_1\leq q_2$, and so each of 
$R(I),R^2(I),R^3(I),R^4(I)$
satisfy the property. If $I=R(J)$, then $q_1\geq q_2$ and so each of
$L(I),L^2(I),L^3(I),L^4(I)$ have the desired property.
The result easily follows.
\hfill{$\Box$}
\medskip

Let $|I|$ denote the length of $I$: if $I=(\frac{p_1}{q_1},\frac{p_2}{q_2})$
then $|I|=1/q_1q_2$.
\begin{cor} If $I$ is a finite Farey interval containing a $4$-aloof number,
then $1/5\leq |L(I)|/|R(I)|\leq 5$.
\label{geom}
\end{cor}
\noindent{\bf Proof.} If $I=(\frac{p_1}{q_1},\frac{p_2}{q_2})$ then
$$\frac{|L(I)|}{|R(I)|}=\frac{|\frac{p_1}{q_1}-\frac{p_1+p_2}{q_1+q_2}|}
{|\frac{p_2}{q_2}-\frac{p_1+p_2}{q_1+q_2}|}=\frac{q_2}{q_1}.$$  \hfill{$\Box$}
\medskip

The following lemma is the key fact which makes the construction work.
\begin{lemma} If $x$ is $4$-aloof and $|\frac kp-x|<\frac 1p$ 
then there is a Farey interval
$(\frac{p_1}{q_1},\frac{p_2}{q_2})$ containing both $x$ and $k/p$ with 
$q_1>c_2\sqrt{p}$
for some universal constant $c_2$.
\label{sqrt}
\end{lemma}
\noindent{\bf Proof.} The Farey intervals nesting down to $x$ decrease geometrically in
size (with scale at most $6$) by Corollary \ref{geom}. 
So there is a Farey interval $I$ 
containing $x$, with  $5/p\leq|I|<30/p$. 
By backing up at most $5$ stages towards the root,
there is a Farey interval $J$  with $I\subset J$ 
such that the distance of $I$ to the endpoints of $J$
is at least $1/p$ (because in the last 5 letters of $w$ 
there is at least one $L$ and one $R$). 
Thus $J$ contains both $x$ and $k/p$. 
Furthermore $\frac 5p\leq|I|\leq |J|\leq \frac{30}{p}\cdot6^5$ again by Corollary 
\ref{geom}.
So if $J=(\frac{p_1}{q_1},\frac{p_2}{q_2})$, we have 
$$\frac{1}{q_1^2}\leq \frac{5}{q_1q_2}=
5|J|\leq \frac{5^26^6}{p},$$
and so taking square roots
$$q_1\geq \frac{\sqrt{p}}{5\cdot 6^3}$$ 
and similarly for $q_2$. \hfill{$\Box$}
\medskip

If the word labelling $J=(\frac{p_1}{q_1},\frac{p_2}{q_2})$ has length $\ell$, then after adding $\ell$ squares to
a $k\times p$ rectangle using the greedy algorithm, we find
the remaining rectangle is $a\times b$, where 
$$\left(\begin{array}{c}a\\b\end{array}\right)=
\left(\begin{array}{cc}p_2&p_1\\q_2&q_1\end{array}\right)^{-1}
\left(\begin{array}{c}k\\p\end{array}\right),$$
and so
$$a=q_1k-p_1p=(\frac{k}{p}-\frac{p_1}{q_1})pq_1\leq |J|pq_1=p/q_2\leq 5\cdot 6^3\sqrt{p},$$
and similarly for $b$.

Now if the aspect ratio $a/b$ or $b/a$ is $x\geq 5$, then backing up one step gives
a rectangle of aspect ratio in $[1,2)$ and sides bounded by $6\cdot 6^3\sqrt{p}$.

\label{scnthick1}
This completes the construction.

\section{Tiling an ``ell''}
\label{ellmethod}
\begin{figure}[htbp]
\centerline{\psfig{figure=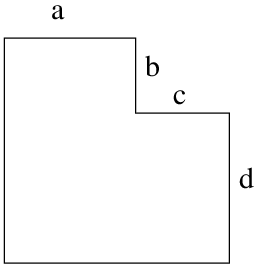,height=1in}}
\caption{\label{el}}
\end{figure}

By an {\bf ell} we mean a rectilinear polygon (polygon with sides parallel to the 
axes) with 6 sides.
We give here a method for tiling (partially) an ell as in Figure \ref{el}
with integer side lengths $a,b,c,d$ as indicated, so that the
remaining untiled portion is an ell with side lengths $a',b',c',d'$ 
of ratios bounded by $8$ (i.e. the ratios are all in $[1/8,8]$).
This construction is a subroutine in the algorithm we will devise in the next section.

Let $L$ be an ell as in Figure \ref{el}.
In what follows we describe an ell with four edge lengths $a,b,c,d$: 
these are the lengths of the four edges
corresponding to the edges marked $a,b,c,d$ of Figure \ref{el}.

Let $a,b,c,d$ be given. If either of $b/a$ or $c/d$ (say $b/a$) 
is $\geq 2$ we add a square
of side $a$ adjacent to the edge of length $a$, 
giving a new ell with $b/a$ reduced by $1$; this does not 
increase the largest ratio of $a,b,c,d$. So in what follows we assume $b/a,c/d<2$.
By symmetry we may assume either $a$ or $b$ is the longest edge.
There are a number of cases to consider\footnote{We apologize for the clumsiness 
of this algorithm.}:

\noindent{\bf Case 1. Suppose $a$ is the longest edge.}

\noindent{\bf Case 1a. Each of the lengths $a,b,c,d$ is in the interval $[d,a]$.} 
We can assume $d<a/3$ or else we are done.

\begin{figure}[htbp]
\centerline{\psfig{figure=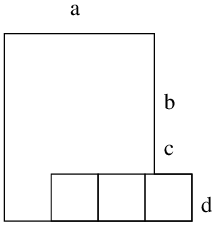,height=1in}}
\caption{\label{da}}
\end{figure}

Add $3$ squares of side length $d$ as in Figure \ref{da}. The new ell has sides
$b,3d-c,d,a+c-3d$. Each of these is in the interval $[d,a]$, since 
$3d-c=d+(2d-c)\geq d$ and $3d-c\leq 3d<a$, and $a>a+(c-3d)=c+(a-3d)>c$.

\noindent{\bf Case 1b. All edges are in $[c,a]$.} 
Suppose also that $c<a/3$. 

\begin{itemize}\item
If $d-c\geq c$, then add a square of side $d$ adjacent to edge $d$. This gives a new
ell with edges $b,d-c,d,a+c-d$, each in $[c,a]$.

\item
If $d-c<c$ and $a+c-2d>c$, add $2$ squares as in Figure \ref{ca}; the remaining ell has
edges $b,2d-c,d,a+c-2d$, and $a>2d>2d-c>d$ by hypothesis 
and so each of these is in $[c,a]$.
\begin{figure}[htbp]
\centerline{\psfig{figure=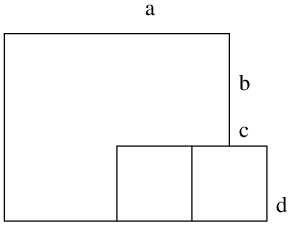,height=1in}}
\caption{\label{ca}}
\end{figure}
\item
If $d-c<c$ and $a+c-2d\leq c$ then $a\leq 2d\leq 4c$ and so $[c,a]\subset[c,4c]$.
\end{itemize}

\noindent{\bf Case 1c. Edges are in $[b,a]$, and $b<a/3$}.

\begin{itemize}\item
If $d-c\geq b$, add a square of size $d$ adjacent to side $d$. The new ell has sides
$b,d-c,d,a+c-d$, each is in $[b,a]$ ($a+c-d=c+(a-d)\geq c\geq b$ and $a+c-d<a$).
\item
If $c>d$, add squares adjacent to $b$ and $d$ as in Figure \ref{ba}; the new ell
has edges $a-b,b,c+b-d,d$, each of which is in $[b,a]$.
(Note $c\geq c+b-d=b+(c-d)\geq b$.)
\begin{figure}[htbp]
\centerline{\psfig{figure=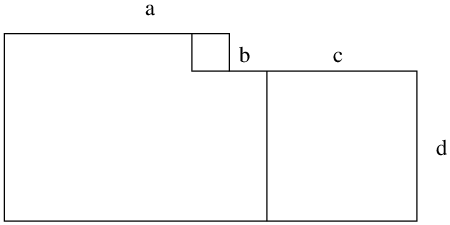,height=1in}}
\caption{\label{ba}}
\end{figure}
\item
If $c\leq d$ and $d-c<b$, add three squares as in Figure \ref{ba2}; the new ell
has edges $a-2b,b,c+2b-d,d$, each of which is in $[b,a]$ (note $c+2b-d<2b<a$).
\begin{figure}[htbp]
\centerline{\psfig{figure=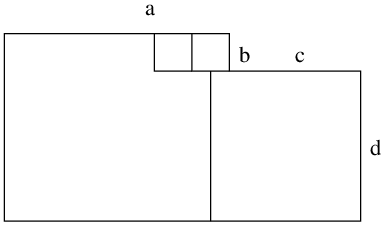,height=1in}}
\caption{\label{ba2}}
\end{figure}
\end{itemize}
\noindent{\bf Case 2. Suppose $b$ is the longest edge.}

\noindent{\bf Case 2a. Edges are in $[a,b]$.} Since $b/a<2$, we are done.

\noindent{\bf Case 2b. Edges are in $[c,b]$.} 
\begin{itemize}\item
If $a\leq 2d$, then
$b<2a\leq 4d<8c$ and so $[c,b]\subset[c,8c]$ and we're done.
\item
If $a> 2d$ and $d-c>c$, then add a square to edge $d$, giving $b,d-c,d,a+c-d$.
\item
If $a>2d$ and $d-c<c$, then add two squares as in Figure \ref{cb}, leaving an ell
with edges $b,2d-c,d,a+c-2d$. Note $2d-c=2(d-c)+c<3c<b$, and $a+c-2d
=a+(c-d)-d<a<b$ so each edge
is in $[c,b]$.
\begin{figure}[htbp]
\centerline{\psfig{figure=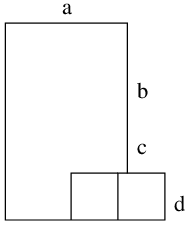,height=1in}}
\caption{\label{cb}}
\end{figure}
\end{itemize}

\noindent{\bf Case 2c. Edges are in $[d,b]$}. 
\begin{itemize}\item
If $a+c-3d\geq d$ add three squares as in Figure \ref{da};
the ell has edges $b,3d-c,d,a+c-3d\in[d,b]$. 
\item
If $a+c-3d<d$ then  $3b/2\leq a+b<a+c<4d$, and so $[d,b]\subset[d,8d/3]$.
\end{itemize}
\medskip

This completes the construction.
Suppose that originally the ratios of $a,b,c,d$ were bounded by $N$. Each square
added at the first step in this algorithm has side length at least the length
of the shortest of $a,b,c,d$. Furthermore, the shortest of $a,b,c,d$ never gets any
shorter. So each square added takes up at least $1/(N^2+N+1)$ of the area.
\label{ellalg}

\section{A better upper bound for rational rectangles.}
\label{scnthick2}
We give in this section a refinement of the construction of 
section \ref{scnthick1}, yielding a logarithmic bound.

The refinement is based on the following theorem, a two-dimensional
version of Theorem \ref{thick}.

Let $C_n$ be the Cantor set of $n$-aloof numbers.
\begin{thm}\label{thick2} For any $M>0$ there is an $n=n(M)$
with the following property.
For any positive real numbers $a,b,c,d$ with ratios bounded by $M$
there exists $t\in [0,a]$ and $r_1,r_2,r_3\in C_n$ such that :
$$r_1=\frac{t}{b+d},~ r_2=\frac{a-t}{b},~r_3=\frac{a+c-t}{d}.$$
\end{thm}

For the proof, see the Section \ref{8}.

\begin{figure}[htbp]
\centerline{\psfig{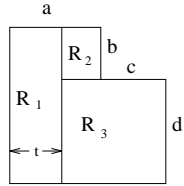}}
\caption{\label{ell}Subdividing an ell into three easy rectangles.}
\end{figure}

The correct interpretation of this theorem is as follows:
Given the ell of Figure \ref{ell}, where the sides $a,b,c,d$ 
have lengths in ratios less than $M$, we can find a $t\in[0,a]$
so that, for the subdivision indicated, the rectangles $R_1,R_2,R_3$
have aspect ratios in $C_{n}$. The quantities $r_i$ of the theorem are the aspect
ratios of the $R_i$ as a function of $t$.

Let $R$ be a $p\times q$ rectangle, $p,q\in\Z$, with $(p,q)=1$ and again $1<q/p\leq2$.
The construction now proceeds as follows.  Let $n=n(M)$ with $M=8$ in Theorem \ref{thick2}.
As before, use Theorem \ref{thick} to divide $R$ into two rectangles $R_1,R_2$,
respectively $p\times k_1$ and $p\times k_2$, 
with $k_i\in\Z$ and $k_i/p$ within $1/p$ of an $n$-aloof number (We assume $n\geq 4$).

Apply the greedy algorithm to $R_1$ and $R_2$ as before, 
until the sides of the untiled rectangles
$R_1',R_2'$ have length $\leq c\sqrt{p}$, and aspect ratios $<2$. (Now the constant $c$
here depends on $n$).
It is easy to arrange that $R_1'$ and $R_2'$ are adjacent, and so
the union of the untiled regions $R_1'\cup R_2'$ then forms an ell.

We claim that we can also arrange so that the ratios of edge lengths of $R_1'$ and $R_2'$
are at most $n$: simply back up the greedy algorithm if necessary for the smaller
of $R_1', R_2'$ until it is approximately the same size as the other. Since the change
in scale between the time the aspect ratio is in $[1,2]$ and the next time
it is in $[1,2]$ 
is at most $n$, this proves the claim.

Using the subroutine of section \ref{ellalg}, we can tile this ell
``easily'' (that is, each square added takes up a definite proportion $1/(n^2+n+1)$ of
the area) until all the ratios of sides are less than $8$.

We then apply Theorem \ref{thick2} with $M=8$: this subdivides the ell  into three
rectangles with aspect ratios in $C_n$. By choosing $t'$ to
be the integer closest to the $t$ of the theorem, we can subdivide the ell into three
rectangles $R_3,R_4,R_5$ each with integer sides of length $\leq c\sqrt{p}$ 
and with
aspect ratios $\frac{t'}{b+d},\frac{a-t'}{b},\frac{a+c-t'}d$ within 
$\frac 1{b+d},\frac 1b,\frac 1d$ respectively of points in $C_n$.

We can now use the greedy algorithm on the $R_3,R_4,R_5$, until the edge lengths
are less than $c^{3/2}p^{1/4}$, and the untiled rectangles $R_3',R_4'$ of $R_3,R_4$
respectively, are adjacent, have aspect ratios $<2$, and are the same size to within a factor of $n$.
(The $R_3',R_4'$ are not necessarily adjacent to $R_5'$).

At the next step the ell formed by $R_3'\cup R_4'$ is tiled using the ell method
of section \ref{ellmethod} and then is subdivided into $3$ ``easy''
rectangles (using Theorem \ref{thick2} again), 
and $R_5'$ is subdivided using Theorem \ref{thick} into two easy rectangles.
As we continue this process, each ell gives rise to an ell and a single rectangle,
and each single rectangle gives rise to an ell.

So the total number of untiled rectangles at the $n$th stage of the construction
is just the $n$th Fibonacci number: letting $f_n,g_n$ 
be the number of ells and single rectangles, after one iteration
we have 
$$\left(\begin{array}{c}f_{n+1}\\g_{n+1}\end{array}\right)=
\left(\begin{array}{cc}1&1\\1&0\end{array}\right)
\left(\begin{array}{c}f_n\\g_n\end{array}\right).$$

As a consequence the number of squares needed to tile is $N$, where:
$$N\leq 2\log_\alpha p+3\log_\alpha (cp^{1/2})+5\log_\alpha (c^{3/2}p^{1/4})+\ldots+
F_k\log_\alpha (c^{2-2^{-k+1}}p^{2^{-k}})+F_{k+1}c',$$
where as before $k\leq \log_2\log_2p$, $c'$ is a bound on the number of squares
needed to tile a rectangle of edge bounded by $2c^2$, and $\alpha$ is a constant
depending on $n$.
This is a convergent geometric series, 
since $F_m\approx \tau^m$ and $\tau\approx1.618<2.$
We have $N\leq C_1\log p$ for some universal constant $C_1$.
This completes Theorem \ref{int}.\hfill{$\Box$}
\label{log}

\section{Proof of Theorem \protect\ref{thick2}}
\label{8}

We will prove a stronger result (Theorem \ref{3}).

Recall that a {\bf gap} of a Cantor set $C\subset\R$ is 
a connected component of $\R-C$.
A Cantor set $C$ is called {\bf $(K,\epsilon)$-thick} if
$C$ is obtained from an interval $J$
by removing successively open subintervals of $J$ which are gaps of $C$, with the 
property: when an gap $I$ is removed from a connected
subinterval $J$, leaving intervals
$I',I''$ on either side with $J=I'\cup I\cup I''$, 
then $|I|\leq \epsilon|J|$ and 
$|I'|/|I''|\in[1/K,K]$.

Recall that $C_n$ is the Cantor set of $n$-aloof numbers.
The following lemma is essentially due to Hall \cite{Hall} (he studied the 
case $n=4$, but his methods extend to any $n$).

\begin{lemma}
For any $\epsilon>0$ there is an integer $n$ such that
$C_n$ is $(3,\epsilon)$-thick.\label{nbig}
\end{lemma}

\begin{thm} Given $M_1,M_2>1$ there exists an $\epsilon>0$ 
with the following property.
Let $E_1,E_2,E_3$ be three $(3,\epsilon)$-thick Cantor sets in 
$\R$ with diameters in ratios bounded by $M_1$.
Let $S$ be the orthogonal projection of $E=E_1\times E_2\times E_3\subset\R^3$ 
to $\R^2$ along
a vector $v\in\R^3$ whose coordinates have ratios in absolute value
bounded by $M_2$.
Then $S$ contains every point in $\R^2$ in the convex hull of $S$ 
which is not within a small neighborhood of the boundary of the convex hull of $S$.
\label{3}
\end{thm}

\noindent{\bf Remark.} Let us show that this theorem implies Theorem \ref{thick2}. 
Take $E_i$ to be the Cantor set $C_n$. For $n$ sufficiently large 
this Cantor set is $(3,\epsilon)$-thick, by Lemma \ref{nbig}.
The set $$\ell=\{(\frac{t}{b+d},\frac{a-t}{b},\frac{a+c-t}{d})~|~t\in[0,a]\}
\subset\R^3$$
is a line segment which passes completely through the convex hull of $E$. The 
direction of $\ell$ is $$v=(\frac{1}{b+d},\frac{-1}{b},\frac{-1}{d}),$$ 
whose coordinate ratios
are bounded in absolute value by $9$ by hypothesis (recall $b/d,d/b\leq 8$). 
We need to show that $\ell$ intersects $E$.

Let $\pi$ be the projection $\pi\colon\R^3\to\R^2$ given by
$$\pi(x,y,z)=(x+\frac{zd}{b+d},y+\frac{zd}{b})$$ 
(which is not the orthogonal projection,
but is orthogonal projection followed by a linear map of bounded
distortion). 
Then $\pi(\ell)$ is the single point $(\frac{a+c}{b+d},\frac{2a+c}{b})$,
which is contained in the square $Q=[1/8,8]\times[3/8,24]$ because the ratios
of any two of $a,b,c,d$ are bounded by $8$, and $Q$ is in turn contained and not close
to the boundary of the convex hull of $\pi(E)$ for $n$ large ($n>24$).
Thus $\pi(\ell)\subset\pi(E)$, and so $\ell$ intersects $E$.
\medskip

\begin{lemma} Fix $N>0$.
If a Cantor set $C\subset \R$ 
is $(3,\epsilon)$-thick for some $\epsilon<1/5$ then for some $k$
we can remove gaps 
$I_1,I_2,\ldots,I_k$ from the convex hull of $C$, leaving subintervals
$J_1,J_2,\ldots,J_{k+1}$, with $|I_i|/|C|<\epsilon$ and for each $i$,
$$\frac{|C|}{5N}<|J_i|\leq \frac{|C|}{N}.$$
\label{subd}
\end{lemma}

\noindent{\bf Proof:} 
Since $C$ is $(3,\epsilon)$-thick with $\epsilon<1/5$, 
each gap $I$ removed from a subinterval $J$ leaves 
two subintervals $I',I''$ of length at least $1/5$th of the length of $J$
(since $|J|=|I'|+|I|+|I''|\leq|I'|+|J|/5+3|I'|$).

So one simply removes gaps from the convex hull
until the remaining subintervals $J_i$ have length between $|C|/N$ and $|C|/5N$.
\hfill{$\Box$}
\medskip

\noindent{\bf Proof of Theorem \ref{3}.}
Our proof remains at a qualitative level for simplicity. In particular we
won't try to estimate the best $\epsilon$.

\begin{figure}[htbp]
\centerline{\psfig{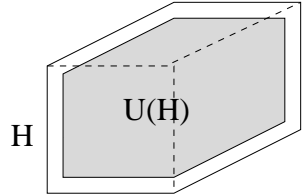}}
\caption{\label{U} The inner neighborhood of a hexagon.}
\end{figure}

Let $I_1,I_2,I_3$ be the convex hulls of $E_1,E_2,E_3$. 
The projection of $I_1\times I_2\times I_3$,
the convex hull of $E$, is a hexagon with opposite sides parallel.
For such a hexagon H, define $U(H)$ to be the set of points in the interior
of $H$ and at distance more than $a/10$ from any boundary edge of length $a$
(see Figure \ref{U}). Call $U(H)$ the {\bf inner neighborhood} of $H$.

Define subdivisions of the $E_i$ as in Lemma \ref{subd} for some large $N$:
so that for $i=1,2,3$ we have
$E_i=E_{i1}\cup E_{i2}\cup\ldots\cup E_{in_i}$ and $|E_{ij}|/|E_i|
\in[1/5N,1/N)$.  The projection of the union of the convex hulls of
the $E_{1i_1}\times E_{2i_2}\times E_{3i_3}$ is a ``stack'' of hexagons
as in Figure \ref{stack}. 
By Lemma \ref{blox} below, $U(E)$
is contained in the union of the inner neighborhoods of the ``blocks''
$E_{1i_1}\times E_{2i_2}\times E_{3i_3}$.
Furthermore each block again satisfies the hypotheses of Lemmas \ref{subd}
and \ref{blox}, and so we can subdivide it again, and repeat. For each point $x\in
U(E)$,
we obtain in this way a sequence of blocks converging to $x$. By compactness $U(E)$
is contained in $\pi(E)$.
\hfill{$\Box$}
\medskip

\begin{lemma} Suppose the same hypotheses as in Theorem 
\ref{thick2}.  For some $N$ large 
define subdivisions $\{E_{1i}\},\{E_{2j}\},\{E_{3k}\},$ 
of $E_1,E_2,E_3$ respectively as in Lemma \ref{subd}.
If $N$ is sufficiently large and $\epsilon$ sufficiently small, 
then $$U(E)\subset \bigcup_{i,j,k}U(E_{1i}\times E_{2j}\times E_{3k}).$$
\label{blox}
\end{lemma}

\noindent{\bf Proof.} The proof by picture is the most illuminating.
The direction of the diagonal edge of the hexagon
in Figure \ref{stack} 
\begin{figure}[htb]
\centerline{\psfig{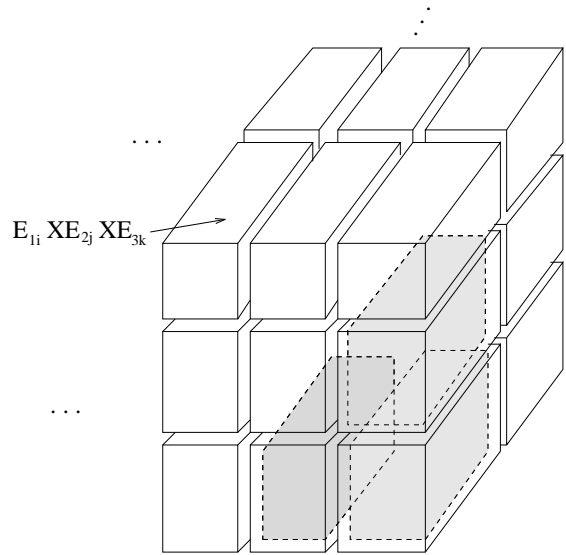}}
\caption{\label{stack} The inner neighborhoods (3 of which are
shaded) of the blocks cover the inner neighborhood of the whole stack.}
\end{figure}
(i.e. the vector $\pi(0,0,1)=(\frac{d}{b+d},\frac db)$) has slope $(b+d)/b$ 
between $1$ and $M_1$.  Using also the fact that the boxes
$E_{1i}\times E_{2j}\times E_{3k}$
have edge-lengths of ratios bounded by a constant ($5M_1$),
their projections are hexagons with edge lengths of ratios bounded by
another constant $C(M_1,M_2)$.

Now if $N$, the approximate
number of blocks per edge in the stack, is sufficiently large and
$\epsilon$ is sufficiently small compared to the size of the smallest
block, we see that the inner neighborhoods of all the
boxes $U(E_{1i}\times E_{2j}\times E_{3k})$ cover all of the convex hull
of $\pi(E)$ except in a small neighborhood of the boundary; in particular they
cover $U(E)$.
\hfill{$\Box$}

\section{Problems}
\begin{prob}
What are the best constants in the upper and lower bound of Theorem \ref{int}?
\end{prob}

Our constructions leave lots of room for improvement in the constant appearing in the
upper bound. The lower bound of $\log_2\max(p,q)$ can also be improved, however.
This is another interesting problem in itself:

\begin{prob}Among all graphs with $m$ edges, which graph $G$ has the largest number
of spanning trees $\kappa(G)$? What is the sup of $\kappa(G)^{1/m}$ over all graphs?
Over planar graphs?
\end{prob}

We used the trivial bound $2^m$ for the number of spanning trees 
of a planar graph. This bound can be improved; N. Young indicated to us an upper bound
of $\lambda^m$ for some $\lambda<2$, which comes from taking into account the
vertex degrees.

On the other hand the $n\times n$ planar grid graph has $\kappa(G)^{1/m}$ converging
to $\approx 1.79$ (see \cite{BP}) as $n\to\infty$, 
and this is the largest value we know of.
So the actual largest value is somewhere in the range $(1.79,2)$.

\begin{prob} How many cubes does it take to tile a $p\times q\times r$
box?
\end{prob}

None of our methods work for this case; even the greedy algorithm is difficult to define.

\section{Appendix}
We give here a proof of Theorem \ref{Per}, which Yuval Peres has
kindly allowed us to include.

Recall the notation: $x\in(0,1)$ and has continued fraction expansion
$x=[0;a_1,a_2,\ldots]$ with $n$th approximants $p_n/q_n$.

Let $S_N(x)$ be the sum of the first $N$ partial quotients of $x$.
Diamond and Vaaler \cite{DV} 
showed that for almost all $x$,
\begin{equation}\label{DV}
S_N(x)=(1+o(1))N\log_2N+\theta\max_{1\leq k\leq N} a_k(x),
\end{equation}
where $\theta\in[0,1]$ (and $\theta$ depends on both $x$ and $N$).

We are interested in $N_\epsilon(x)=\min\{N : |q_Nx-p_N|<\epsilon\}$.
By a result of Khinchin and Levy (cf \cite{Bill}), for almost every $x$
$$\frac1N\log|q_Nx-p_N|\to -c_1= -\frac{\pi^2}{12\log 2}.$$
By discarding a set of measure $\delta$ for any small $\delta>0$,
this convergence is uniform, i.e. on a set $A\subset(0,1)$
with $\mu(A)>1-\delta$, we have 
$$\sup_A\{\frac 1N\log |q_Nx-p_N|+c_1\}\to 0.$$ 
We conclude that 
$\log\frac 1\epsilon/N_\epsilon(x)$ converges uniformly to $c_1$
on $A$.

Using the Gauss-Kuz'min measure we have 
$\mu\{x : a_i(x)\geq k\}< c_2/k$ for a constant $c_2$;
and so on $A$ we have (using $M_\epsilon=\frac{1}{c_1}\log\frac{1}{\epsilon}$):
\begin{eqnarray*}\mu\left\{x\in A \mid \max_{1\leq k\leq M_\epsilon} a_k \geq M_\epsilon\sqrt{\log 
M_\epsilon}\right\}&\leq& M_\epsilon\mu\left\{x\in A \mid a_1\geq M_\epsilon\sqrt{\log
M_\epsilon}\right\}\\&\leq&c_2/\sqrt{\log M_\epsilon}.
\end{eqnarray*}
Letting $B$ be the complement of this set in $A$, using (\ref{DV}) and $M_\epsilon=
\frac1{c_1}\log\frac 1\epsilon$ we have for all $x\in B$
\begin{eqnarray*}
T_\epsilon(x)=S_{N_\epsilon(x)}&=&(1+o(1))N_\epsilon(x)\log N_\epsilon(x)+\theta
M_\epsilon\sqrt{\log M_\epsilon}\\
&=&(1+o(1))c_3\log\frac{1}{\epsilon}\log\log\frac{1}{\epsilon}\\
\end{eqnarray*}
for some constant $c_3>0$.
\hfill{$\Box$}

\end{document}